\documentclass[a4paper]{amsart}
\usepackage{latexsym}
\usepackage{graphics}
\newtheorem{theorem}{Theorem}
\newtheorem*{corollary}{Corollary}
\newtheorem*{lemma}{Lemma}
\newtheorem*{proposition}{Proposition}
\theoremstyle{definition}
\newtheorem*{definition}{Definition}
\theoremstyle{remark}
\newtheorem*{remark}{Remark}
\numberwithin{equation}{section}
\newcommand{\R}{\mathbb{R}}
\newcommand{\C}{\mathbb{C}}
\newcommand{\calO}{{\mathcal O}}    
\newcommand{\calG}{{\mathcal G}}    
\newcommand{\Aut}{\mathrm{Aut}}
\newcommand{\Out}{\mathrm{Out}}
\newcommand{\im}{{\mathrm{im}}}                
\newcommand{\orbchi}{\chi^{\mathrm{orb}}}      
\newcommand{\nohyphen}{\discretionary{}{}{}}   
\newcommand{\denom}{\frac{1}{\sqrt{2\pi}}}
\begin{document}
\title{The Euler characteristic of graph complexes via Feynman diagrams}
%
\author{Ferenc Gerlits 
}                     
%
%
\address{Alfr\'ed R\'enyi Institute of Mathematics, Re\'altanoda
u.\ 13--15, Budapest, Hungary}
%
\date{Dec.~4, 2004}
%
\begin{abstract}
We prove several claims made by Kontsevich about the orbifold
Euler characteristic of the three types of graph homology
introduced by him.  For this purpose, first we develop a
simplified version of the Feynman diagram method, which requires
integrals in one variable only, to obtain the generating functions
as asymptotic expansions of certain Gaussian integrals.  Finally,
following Penner, we relate these integrals to the gamma function
in order to compute the individual coefficients of the generating
functions.
\end{abstract}
\maketitle

\section{Introduction} \label{sec:intro}
In two papers (\cite{Konts-Formal}, \cite{Konts-Feynman})
published in 1993 and 1994, Kontsevich constructed a family of
objects called \emph{graph complexes}.  These are chain complexes
of vector spaces, where each vector space is spanned by a set of
finite graphs equipped with an orientation.  He computed the
orbifold Euler characteristic of these chain complexes, and stated
the answers; however, he did not give any proofs except for a one
sentence reference to Feynman diagrams.

In Sections~\ref{sec:eulerchar}, \ref{sec:graphhom}
and~\ref{sec:integrals} we recall the definitions of the graph
complexes and their orbifold Euler characteristic, and prove most
of Kontsevich's claims about them. But first, we introduce the two
technical tools used in the proof. In Section~\ref{sec:graphs}, we
develop a simplified version of the \emph{Feynman diagram\/}
apparatus (also called the method of \emph{Gaussian integrals\/}).
This method, well known to physicists but mostly unknown to
mathematicians, can be used to construct generating functions for
various problems involving the counting of graphs. A good
exposition of the general method can be found in \cite{BIZ}. Our
version is based on \cite{Penner}.  It is limited in that it cannot
take into account the number of boundary components of ribbon graphs
(referred to as the ``topological expansion'' in the physics literature),
but it is much simpler as it requires integration over $\R$ only,
instead of over a space of matrices.  In Section~\ref{sec:analysis},
we recall some basic facts about \emph{asymptotic expansions}, which are
used to extract information about the individual terms of the
generating functions obtained by the Feynman diagram method.

The emphasis is on proving the claims made by Kontsevich about the
orbifold Euler characteristic of the commutative, associative and
Lie graph complexes.  However, our simplified Feynman diagram
method can be applied to graph complexes based on other cyclic
operads (see \cite{KarenJim} or \cite{Swapneel}), as well as to
other problems involving the counting of graphs.

\section{Acknowledgements}
This work has grown out of a seminar organized by Karen Vogtmann
in Fall 2000 at Cornell University, with the goal of understanding
Kontsevich's graph homology.  It is based on Chapter~5 of the
author's Ph.D.\nocite{mythesis} dissertation, which could not have
been written without Swapneel Mahajan's help.

\section{Enumeration of graphs}\label{sec:graphs}
As our first example, we consider the following counting problem.
Let $T_{2e}$ be the set of graphs with $e$ edges such that each
vertex has valence at least 3.  (The \emph{valence\/} of a vertex
is the number of edges incident to it.)  For each $G\in T_{2e}$,
let $\Aut(G)$ be the set of automorphisms of $G$.  The goal is to
evaluate the sum
\begin{equation}
  \sum_{G\in T_{2e}} \frac{1}{|\Aut(G)|}.
\end{equation}

\begin{figure}[h]
\resizebox{0.9\textwidth}{!}{%
  \includegraphics{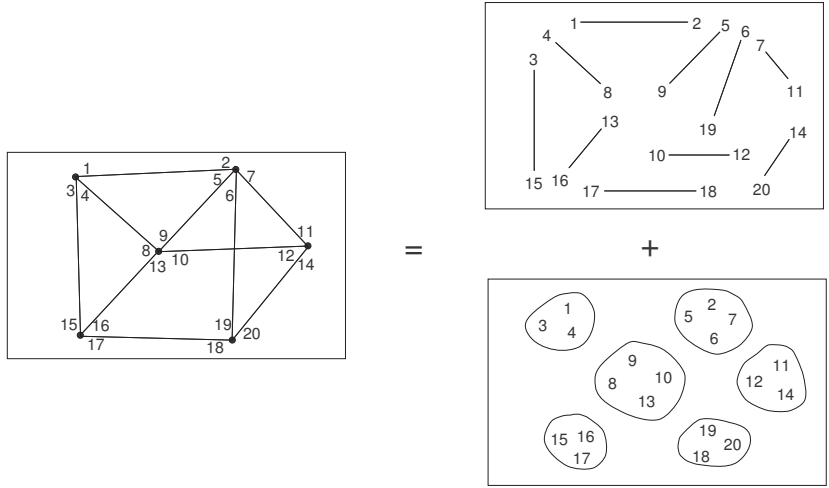}
}
\caption{Labeled graph = chord diagram + partition}
\label{fig:ExampleGraph}
\end{figure}

The first observation is that we are, in fact, counting the number
of \emph{labeled graphs.\/}  We will say that a graph with $e$
edges is \emph{labeled} if the $2e$ half-edges are numbered from
$1$ through $2e$ (see Figure~\ref{fig:ExampleGraph}).  Fix a graph
$G$, and consider all its labelings, i.e., the (isomorphism
classes of) labeled graphs with $G$ as the underlying unlabeled
graph.  The symmetric group $S_{2e}$ acts on these labelings, and
the size of the stabilizer of any given labeling is equal to the
number of automorphisms of $G$. Thus if $T^{\ell}_{2e}$ denotes
the set of labeled graphs with $e$ edges such that each vertex is
at least trivalent, then
\begin{equation} \label{eq:orbit-stabilizer}
  \sum_{G\in T_{2e}} \frac{1}{|\Aut(G)|} =
  \sum_{G\in T_{2e}} \frac{|\{\text{labelings of $G$}\}|}{(2e)!} =
  \sum_{G\in T^{\ell}_{2e}} \frac{1}{(2e)!},
\end{equation}
i.e., the sum to be computed is the number of labeled graphs
$|T^{\ell}_{2e}|$ divided by $(2e)!$.

Second, notice that we can break up a labeled graph into two
separate objects: a pairing of the numbers $1$, $2$, \ldots, $2e$,
which we will call a \emph{chord diagram}, describes the edges;
and a partition of these $2e$ numbers into groups of three or more
tells us what the vertices are.  Conversely, any pair of a chord
diagram and such a partition corresponds to a labeled graph.
Therefore, if $Ch_n$ denotes the number of chord diagrams on $n$
numbers and $Par^{3+}_n$ the number of set partitions of $n$
distinct objects into subsets of size three or more, then
\begin{equation}
  \label{eq:ChTimesPar}
  |T^{\ell}_{2e}| = Ch_{2e} Par^{3+}_{2e}.
\end{equation}

Each of the quantities on the right hand side are fairly easy to
compute. In the case of chord diagrams, this could be done
directly, but in both cases the simplest solution is to find their
exponential generating functions $Ch(x) = \sum_{n=0}^{\infty} Ch_n
x^n /n!$ and $Par^{3+}(x) = \sum_{n=0}^{\infty} Par^{3+}_n x^n
/n!$.

Namely, let $Q$ denote some class of objects; following
\cite{Joyal}, we will call it a \emph{species}.  Let $Q(x) =
\sum_{n=0}^{\infty} Q_n x^n /n!$ be the corresponding exponential
generating function, where $Q_n$ is the number of objects of type
$Q$ on $n$ numbers, vertices etc.  Consider the \emph{exponential
species} $\exp(Q)$, the objects of which are the disjoint unions
of any number of objects of type $Q$.  Then the exponential
generating function $(\exp(Q))(x)$ is equal to the formal
exponential
\begin{equation}
  \exp(Q(x)) := 1 + Q(x) + \frac{1}{2}(Q(x))^2 + \cdots
                  + \frac{1}{n!}(Q(x))^n + \cdots
\end{equation}
of the power series $Q(x)$, provided that any object of type
$\exp(Q)$ can be uniquely decomposed into its sub-objects of type
$Q$. In particular, we must have $Q_0 = 0$, so the sum above is
finite in each coefficient. This is an application of Joyal's
theory of combinatorial species; details can be found in
\cite{redbook} or \cite{vanLint}.

Now, let $E^{2}$ be the species of two-element sets, i.e.,
$E^{2}_2 = 1$ and $E^{2}_n = 0$ for all $n\ne 2$.  The exponential
generating function is $E^{2}(x) = x^2/2$.  On the other hand, we
have $Ch = \exp(E^{2})$, and thus
\begin{equation}
  Ch(x) = e^{x^2/2}.
\end{equation}
Similarly, if $E^{3+}$ denotes the species of sets of size greater
than or equal to three, then we have
\begin{equation}
  E^{3+}(x) = \frac{x^3}{3!} + \frac{x^4}{4!} + \cdots
       = e^x - 1 - x - x^2/2.
\end{equation}
Since $Par^{3+} = \exp(E^{3+})$, we immediately obtain
\begin{equation}
  Par^{3+} (x) = e^{e^x - 1 - x - x^2/2}.
\end{equation}

Finally, we would like to find the exponential generating function
of $|T^{\ell}|$, given those of $Ch$ and $Par^{3+}$, using
equation~(\ref{eq:ChTimesPar}).  This is called the
\emph{Cartesian product} of the two species, denoted by
$|T^{\ell}| = Ch \times Par^{3+}$.  There does not seem to be any
general method for computing the generating function of a
Cartesian product; however, in this particular case, when one of
the factors is $Ch$, the following lemma solves the problem.
Readers familiar with the the theory of Feynman diagrams will
recognize it as a one-dimensional version of the calculation
underlying Wick's Lemma.

\begin{lemma} \label{lem:Wick's}
$\displaystyle Ch(x) = \frac{1}{\sqrt{2\pi}}
\int_{-\infty}^{\infty} e^{-y^2/2} e^{xy}\,dy$.
\end{lemma}

\begin{proof}
Straightforward.
\end{proof}

Another way of stating the lemma is
\begin{equation}
  Ch_n =
  \frac{1}{\sqrt{2\pi}} \int_{-\infty}^{\infty} e^{-y^2/2}y^n\,dy
        \mathrm{\qquad for\ all\ } n.
\end{equation}

Note that this formula holds for odd values of $n$, as well: in
this case, both sides are zero.  If we could exchange the order of
the summation and the integration, we would have

\begin{equation}
\label{eq:Kontsevich}
  |T^\ell|(x) = \frac{1}{\sqrt{2\pi}} \int_{-\infty}^{\infty}
              e^{-y^2/2}e^{e^{xy} - 1 - xy - x^2y^2/2}\,dy.
\end{equation}

Although this formula needs explanation (the series $|T^\ell|(x)$
does not converge: its radius of convergence is zero), it is
indeed valid in a limited sense, as we shall see below.

\section{Asymptotic expansions}\label{sec:analysis}
Power series with zero radius of convergence can be treated
analytically by considering them as the asymptotic expansion of
some function.  We briefly recall the definition and basic
properties of asymptotic expansions; details can be found in
\cite{deBruijn}, for example.

\begin{definition}
Let $f$ be a function defined on a subset of $\R$ including a
(possibly one-sided) neighborhood of~0, with values in $\R$ or
$\C$. We say that the formal power series
$\sum_{n=0}^\infty{a_n}{x^n}$ is an \emph{asymptotic expansion} of
$f$ around $x=0$ and we write
\begin{equation}
  f(x) \approx \sum_{n=0}^\infty a_n x^n \qquad\text{as $x\to 0$}
\end{equation}
if the following condition is satisfied for each $N\ge 1$:
\begin{equation} \label{eq:AsympDef}
  f(x) = \sum_{n=0}^{N-1} a_n x^n + \calO(x^{N})
             \qquad\text{as $x\to 0$.}
\end{equation}
This means that for each fixed $N$, the quotient $(f(x) -
\sum_{n=0}^{N-1} a_n x^n)/x^N$ is bounded in some neighborhood
of~$0$.
\end{definition}

Note that the neighborhoods we take in (\ref{eq:AsympDef}) may
shrink down to zero as $N\to \infty$, in which case our definition
does not say anything about the tail of the series
$\sum_{n=0}^\infty a_n x^n$ for any fixed non-zero $x$.  It is
possible that $f(x) \approx \sum_{n=0}^\infty a_n x^n$, but the
series does not converge for any non-zero value of $x$.  Even when
the series does converge, its sum is not necessarily equal to
$f(x)$; for example, $e^{-1/x} \approx 0 + 0x + 0x^2 + \cdots$ as
$x\to 0^+$.  Also note that the function $f$ does not need to be
differentiable, or even continuous.

On the other hand, if $f$ is a $C^\infty$ function, and its Taylor
series has a positive radius of convergence, then the Taylor
series is an asymptotic expansion of $f$. A function may not have
an asymptotic expansion, but if it does, it is unique: one can
show by induction on (\ref{eq:AsympDef}) that the two series are
equal term by term. All the usual operations on formal power
series carry over to asymptotic expansions: given the asymptotic
expansions of $f$ and $g$, the functions $f+g$, $fg$, $e^f$ etc.\
all have asymptotic expansions, which are equal to the sum,
product, exponential etc.\ of those of $f$ and $g$, as formal
power series.

If we consider it as an asymptotic expansion, any power series can
be exchanged with the integral in the Lemma of the previous
section. More precisely, the following is true.

\begin{proposition} \label{prop:Exchange}
Suppose that $h\colon \R \to \C$ has the asymptotic expansion
around zero $h(x) \approx \sum_{n=0}^\infty{a_n}{x^n}$, and the
integral $f(x) = \frac{1}{\sqrt{2\pi}} \int_{-\infty}^\infty
e^{-y^2/2} h(xy)\, dy$ converges for all $x$ in some neighborhood
of zero. Then $f$ has the following asymptotic expansion around
zero:
\begin{equation}
\begin{aligned}
  f(x) &\approx \sum_{n=0}^\infty \left(
   \frac{1}{\sqrt{2\pi}} \int_{-\infty}^\infty e^{-y^2/2} y^n\,dy
   \right) a_n x^n \qquad\text{as $x\to 0$} \\
   &\approx \sum_{n=0}^{\infty} Ch_n a_n x^n \qquad\text{as $x\to 0$.}
\end{aligned}
\end{equation}
\end{proposition}

\begin{proof}
Since $h(x) = \sum_{n=0}^{N-1}{a_n}{x^n} + \calO(x^N)$, we have
\begin{align}
 \frac{1}{\sqrt{2\pi}} &\int_{-\infty}^\infty e^{-y^2/2} h(xy)\,dy \\
 &= \frac{1}{\sqrt{2\pi}} \int_{-\infty}^\infty e^{-y^2/2}
    \biggl( \sum_{n=0}^{N-1}{a_n}{x^n}y^n + \calO(x^Ny^N) \biggr)\,dy \\
 &= \frac{1}{\sqrt{2\pi}} \int_{-\infty}^\infty e^{-y^2/2}
       \sum_{n=0}^{N-1}{a_n}{x^n}y^n \, dy +
    \frac{1}{\sqrt{2\pi}} \int_{-\infty}^\infty e^{-y^2/2}
       \calO(x^Ny^N)\, dy \\
 &= \sum_{n=0}^{N-1} \left(
    \frac{1}{\sqrt{2\pi}} \int_{-\infty}^\infty e^{-y^2/2} y^n\,dy
    \right) a_n x^n + \calO(x^N),
\end{align}
because $\int_{-\infty}^\infty e^{-y^2/2} y^N\,dy$ is finite for
any fixed value of $N$.
\end{proof}

Or in other words:

\begin{theorem} \label{thm:Exchange}
Let $P$ be any species.  If there is a function $h\colon\R\to\C$
such that $P(x)$ is the asymptotic expansion of $h$ around zero,
then
\begin{equation}
  \frac{1}{\sqrt{2\pi}} \int_{-\infty}^\infty
     e^{-y^2/2} h(xy)\, dy
  \approx (Ch \times P)(x) \qquad \text{as $x\to 0$}
\end{equation}
as long as the the integral converges for all $x$ in some
neighborhood of zero.
\end{theorem}

We will also need the following more general version, which can be
proved in the same way.

\begin{corollary} \label{cor:DoubleExchange}
Let $h\colon\R\times\R \to \C$ be defined on $U\times\R$, where
$U$ is some neighborhood of zero, and suppose there is a sequence
$p_n$ of polynomials such that
\begin{equation} \label{eq:PolynomialAsympEq}
  h(x,y) = \sum_{n=0}^{N-1} p_n(y)x^n + \calO(x^Ny^{M(N)})
     \qquad \text{as $x\to 0$}
\end{equation}
holds for each $N\ge 1$, where $M(N)$ is an integer depending only
on $N$.  Suppose further that the integral $f(x) =
\frac{1}{\sqrt{2\pi}} \int_{-\infty}^\infty e^{-y^2/2} h(x,y)\,
dy$ converges for all $x\in U$.  Then $f$ has an asymptotic
expansion around zero, which can be obtained from
(\ref{eq:PolynomialAsympEq}) by replacing each occurrence
of\/~$y^k$ by the number $Ch_k$.
\end{corollary}

\section{The orbifold Euler characteristic}\label{sec:eulerchar}
Our goal is the compute the orbifold Euler characteristic of
certain graph complexes.  A \emph{graph complex} is a chain
complex
\newcommand{\arrow}{\longrightarrow}
\[
  \cdots \arrow C_{k+1} \overset{\partial_{k+1}}{\arrow}
         C_k \overset{\partial_k}{\arrow} C_{k-1} \arrow \cdots
         \arrow C_1 \overset{\partial_1}{\arrow} C_0 \arrow 0
\]
where each $C_k$ is a vector space (in this paper, always over
$\R$) spanned by some set of graphs.  In our examples, graphs in
$C_k$ have $k$ vertices.  Given a graph complex $C_*$ of finite
length where each $C_k$ is spanned by a basis corresponding to the
finite list of graphs $G_{k1}$, $G_{k2}$, $\ldots$, $G_{kn_k}$,
its \emph{orbifold Euler characteristic} is defined as
\begin{equation}
  \orbchi(C_*) = \sum_{k\ge 0} (-1)^k \sum_{l=1}^{n_k}
                   \frac{1}{|\Aut(G_{kl})|}.
\end{equation}
Or, if $\calG$ is the set of all the graphs $G_{kl}$ and $v(G)$ is
the number of vertices of the graph $G$, then we can write
\begin{equation}
  \orbchi(C_*) = \sum_{G\in\calG} \frac{(-1)^{v(G)}}{|\Aut(G)|}.
\end{equation}

Let $C^{(n)}_*$ be the graph complex where $C^{(n)}_k$ is spanned
by all graphs of valence 3 or more on $k$ vertices such that the
number of edges minus the number of vertices is $n-1$.  (If the
graph is connected, this means that there are $n$ independent
loops in the graph, i.e., the fundamental group of the graph has
rank~$n$.)  We want to compute the generating function
\begin{equation}
  \orbchi(C_*)(t) = \sum_{n\ge 2} \orbchi(C^{(n)}_*) t^n
\end{equation}
of the orbifold Euler characteristic of this graph complex.

First, consider the power series
\begin{equation}
  h(t,x,y) = \exp(-\frac1t E^{3+}(xy)),
\end{equation}
and replace each occurrence of $y^k$ by the number $Ch_k$ (note
that this eliminates all the odd powers of $x$).  We get a power
series in $t^{-1}$ and $x$, which is almost the same as the series
$|T^{\ell}|(x) = (Ch \times Par^{3+})(x)$ considered in
Section~\ref{sec:graphs}, but with an extra $-t^{-1}$ for each
vertex of the graphs.  More precisely, the coefficient of $x^{2e}
t^{-v}$ is the sum $(-1)^v \sum 1/|\Aut(G)|$ over graphs with $e$
edges and $v$ vertices.

Now, substitute $x = t^{1/2}$ to obtain the power series
\begin{equation} \label{eq:DoublePowerSeries}
  h(t,y) = \exp(-\frac1t E^{3+}(t^{1/2}y)).
\end{equation}
If we substitute $Ch_k$ for each occurrence of $y^k$, we get a
power series in $t$, where the coefficient of $t^n$ is the sum
$\sum (-1)^{v(G)}/|\Aut(G)|$ over graphs with $n =
\#\{\text{edges}\} - \#\{\text{vertices}\}$.  Note that $E^{3+}_0
= E^{3+}_1 = E^{3+}_2 = 0$ ensures that all the exponents of $t$
are non-negative, and that if we order $h(t,y)$ by exponents of
$t$, then the coefficient of $t^n$ is a polynomial in $y$, with
exponents between $2n$ and $6(n-1)$.

So we have proved the following: if we replace each occurrence
of\/ $y^k$ by $Ch_k$ in the power series $h(t,y)$ defined above in
(\ref{eq:DoublePowerSeries}), we get $t^{-1}$ times the generating
function $\orbchi(C_*)(t)$.

Now it is a simple matter to check that the conditions of the
Corollary to Theorem~\ref{thm:Exchange} are satisfied, and hence
\begin{equation}
  \denom\int_{-\infty}^\infty e^{-y^2/2}h(t,y)\,dy
    \approx t^{-1} \orbchi(C_*)(t) \qquad \text{as $t\to 0$}.
\end{equation}

The same method works for many different types of graph complexes.

\begin{theorem} \label{thm:Kontsevich}
Consider the graph complex $C_*$ consisting of graphs with some
structure at each vertex, represented by the species Q, where
$C_k$ is spanned by all such $Q$-graphs on $k$ vertices.  If\/
$Q_0 = Q_1 = Q_2 = 0$, then the orbifold Euler
char\-ac\-ter\-is\-tic of the graph complex can be obtained from
the power series
\begin{equation} \label{eq:GeneralDoublePowerSeries}
  h(t,y) = \exp(-\frac1t Q(t^{1/2}y))
\end{equation}
by replacing each occurrence of $y^k$ by $Ch_k$.

If, moreover, there is a function (also denoted by $h(t,y$)) which
has the asymptotic expansion (\ref{eq:GeneralDoublePowerSeries})
and which satisfies the conditions of the Corollary to
Theorem~\ref{thm:Exchange}, then the orbifold Euler characteristic
is an asymptotic expansion:
\begin{equation}  \label{eq:GeneralKontsevich}
   \denom \int_{-\infty}^{\infty} e^{-y^2/2} h(t,y)\,dy
   \approx t^{-1}\orbchi(C_*)(t) \qquad \text{as $t \to 0.$}
\end{equation}
\end{theorem}

\section{The three graph homologies}\label{sec:graphhom}
In \cite{Konts-Formal}, Kontsevich constructed three graph
complexes, spanned by what we will call commutative, associative
and Lie graphs.  The homology of the chain complex, i.e., the
quotient $\ker(\partial) / \im(\partial)$ is Kontsevich's
\emph{graph homology}.

\emph{Commutative graphs\/} are the simplest: they are finite
graphs of valence 3 or higher, i.e., of the type considered so
far.  The species $Q$ at each vertex is $E^{3+}$.  Classes in
commutative graph homology play a role in the theory of finite
type invariants of homology 3-spheres \cite{BGRT}.

\emph{Associative graphs\/} have the additional structure of a
cyclic ordering of the half-edges at each vertex.  These are the
same as the \emph{ribbon graphs} or \emph{fat graphs} of
Teichm\"uller theory.  Associative graph homology is the direct
sum of the homology of certain moduli spaces of mapping class
groups \cite{KarenJim}. The species at each vertex is $C$, the
species of cyclic orderings; its terms are $C_n = (n-1)!$.

\emph{Lie graphs\/} have a planar trivalent tree at each vertex,
subject to some relations corresponding to the anti-symmetry and
Jacobi identity of Lie algebras.  The corresponding graph homology
computes the homology (with trivial real coefficients) of the
group $\Out(F_n)$.  The details can be found in \cite{KarenJim};
for now we only need the number (more precisely: the dimension of
the quotient vector space) of such trees with $n$ leaves, which is
$Lie_n = (n-2)!$.

Other types of graph homologies can also be considered, and
Kontsevich's machinery, including the Euler characteristic
computations explained here, can be applied to them; but only
these three are known so far to have applications in other parts
of mathematics.

There are two complications.  First, Kontsevich considers
\emph{connected} graphs only, thereby obtaining the primitive
homology of the corresponding topological object.  This is easy to
remedy: the logarithm of the generating function for the orbifold
Euler characteristic of all graphs yields that of connected graphs
only.

The second problem is that Kontsevich's graphs come with a sense
of \emph{orientation,} and the chain groups $C_*$ are spanned by
those graphs only which do not have any orientation-reversing
automorphisms.  By a lucky coincidence, in the associative and Lie
cases, removing the graphs with orientation-reversing
automorphisms does not change the homology or the orbifold Euler
characteristic; this is shown in \cite{KarenJim}.  In the
commutative case, unfortunately, it does: already in the case of 3
independent loops, we get $\orbchi = 0$ if we include graphs with
orientation-reversing automorphisms, and $\orbchi = -1/48$ if we
do not.  According to Kontsevich, we do get an approximately
correct answer, though, because the typical graph does not have
any automorphisms.

\newcommand{\dis}{\displaystyle}
\renewcommand{\arraystretch}{2.25}
\begin{table}
\caption{The orbifold Euler characteristic of Kontsevich's
                 graph complexes}
\label{table:OrbChi}
\begin{tabular}{c|c|c|c}
 & \multicolumn{3}{c}{orbifold Euler characteristic} \\
\# loops & $\orbchi_c$ (commutative) & $\orbchi_a$ (associative) & $\orbchi_l$ (Lie) \\
\hline 2 & $\dis {\frac{1}{12}}$ & $\dis {\frac{1}{12}}$ &
$\dis {-\frac{1}{24}}$ \\

3 & $\dis {0}$ & $\dis {0}$ &
$\dis {-\frac{1}{48}}$ \\

4 & $\dis {-\frac{1}{360}}$ & $\dis {-\frac{1}{360}}$ &
$\dis {-\frac{161}{5760}}$ \\

5 & $\dis {0}$ & $\dis {0}$ &
$\dis {-\frac{367}{5760}}$ \\

6 & $\dis {\frac{1}{1260}}$ & $\dis {\frac{1}{1260}}$ &
$\dis {-\frac{120257}{580608}}$ \\

7 & $\dis {0}$ & $\dis {0}$ &
$\dis {-\frac{39793}{45360}}$ \\

8 & $\dis {-\frac{1}{1680}}$ & $\dis {-\frac{1}{1680}}$ &
$\dis {-\frac{6389072441}{1393459200}}$ \\

9 & $\dis {0}$ & $\dis {0}$ &
$\dis {-\frac{993607187}{34836480}}$ \\

10 & $\dis {\frac{1}{1188}}$ & $\dis {\frac{1}{1188}}$ &
$\dis {-\frac{5048071877071}{24524881920}}$ \\

11 & $\dis {0}$ & $\dis {0}$ &
$\dis {-\frac{9718190078959}{5748019200}}$ \\
\end{tabular}
\end{table}
Table~\ref{table:OrbChi} shows the orbifold Euler characteristics
$\orbchi_c$, $\orbchi_a$ and $\orbchi_l$ of Kontsevich's
commutative, associative and Lie graph complexes, respectively,
for small numbers of independent loops.  The results are after
taking the logarithm, but without removing graphs with
orientation-reversing automorphisms. It was computed using
(\ref{eq:GeneralDoublePowerSeries}) and the computer algebra
software Mathematica.

A striking feature of this ``experimental data'' is the equality
of the commutative and the associative orbifold Euler
characteristics.  Furthermore, if we type these numbers into the
On-Line Encyclopedia of Integer Sequences \cite{sequenceDB}, we
find that
\begin{equation}\label{eq:Bernoulli}
  (\orbchi_c)_n = (\orbchi_a)_n = \frac{B_{n}}{n(n-1)},
\end{equation}
where the $B_n$ are the Bernoulli numbers.  By using asymptotic
expansions, we can prove these facts.

\begin{remark}
The numbers in Table~\ref{table:OrbChi} are not new. The
associative case (including (\ref{eq:Bernoulli}) for
$\orbchi_a\/$) follows from Harer and Zagier's
computations~\cite{HarZag} of the Euler characteristic of mapping
class groups; the power series $\orbchi_l(t)$ was computed by
Smillie and Vogtmann~\cite{KarenChi} as the generating function of
the orbifold Euler characteristic of $\Out(F_n)$.  Our paper
confirms their computations by describing another, simpler,
derivation of these results.
\end{remark}

\section{Evaluating the integrals}\label{sec:integrals}
The following computation in the commutative case is the same as
Penner's \cite{Penner}.  The other two sections are (unfortunately
incomplete) attempts at extending his method to the associative
and Lie graph complexes, as well.
\subsection{The commutative case}
The commutative series $Q(x) = E^{3+}(x)$ is the Taylor expansion
of $f(x) = e^x - 1 - x - x^2/2$ around~$0$.  We can apply
Theorem~\ref{thm:Kontsevich} directly with $h(t,y) = \exp(-\frac1t
f(t^{1/2}y))$, which gives us

\begin{equation}
\begin{aligned}
   t^{-1}\orbchi_c(t) &\approx \log\biggl(
       \frac{1}{\sqrt{2\pi}} \int_{-\infty}^{\infty}
       e^{-y^2/2}\exp\left(-\frac{1}{t}\bigl(
       e^{y\sqrt{t}} - 1 - y\sqrt{t} - y^2t/2\bigr)\right)
       \,dy \biggr) \\
   &= \log\biggl(
      \frac{1}{\sqrt{2\pi t}} \int_{-\infty}^{\infty}
          e^{-(e^u - 1 - u)/t} \,du \biggr) \\
   &= \log\left(\frac{1}{\sqrt{2\pi t}} \int_{-\infty}^{\infty}
          e^{-\frac{e^u}{t}} e^{\frac{1}{t}} e^{\frac{u}{t}} \,dy
          \right).
\end{aligned}
\end{equation}
Here we substituted $u = y\sqrt{t}$.  Now substitute $z = e^u/t$
to obtain
\begin{equation}\label{eq:CommOrbChi2}
\begin{aligned}
   t^{-1}\orbchi_c(t) &\approx \log\left(
   \frac{e^{1/t}}{\sqrt{2\pi t}} \int_{-\infty}^{\infty}
        e^{-z} (tz)^{\frac{1}{t}} \frac{dz}{z} \right) \\
   &= \log\left(
        \frac{(et)^{1/t}}{\sqrt{2\pi t}} \int_{-\infty}^{\infty}
        e^{-z} z^{\frac{1}{t} - 1} \,dz \right) \\
   &= \log\left( \frac{(et)^{1/t}}{\sqrt{2\pi t}}
                  \,\Gamma\left(\frac{1}{t}\right) \right),
\end{aligned}
\end{equation}
where in the last step we applied the integral formula $\Gamma(x)
= \int_{-\infty}^\infty e^{-z} z^{x-1} \,dz$ which is one of
several equivalent definitions of the gamma function. It is well
known that $\Gamma(z)$ has the asymptotic expansion around
infinity given by Stirling's formula
\begin{equation}
 \Gamma(z) \approx \sqrt{\frac{2\pi}{z}} (z/e)^z e^{J(z)}
 \qquad\text{as $z\to \infty$,}
\end{equation}
(see, for example, \cite{Ahlfors}), where
\begin{equation}
  J(z)=\sum_{n=1}^{\infty} \frac{B_{2n}}{2n(2n-1)}
                            \left(\frac{1}{z}\right)^{2n-1}
\end{equation}
and the $B_{2n}$ are the Bernoulli numbers. If we apply this
formula to $z = 1/t$, take the logarithm, and compare terms, we
get
\begin{equation} \label{eq:CommOrbChi}
 \chi^{\text{orb}}_c(t) = \sum_{n=1}^{\infty}
                      \frac{B_{2n}}{2n(2n-1)} \, t^{2n}.
\end{equation}
In other words, the orbifold Euler characteristic of $C^{(m)}_*$
is zero for $m$ odd, and it is given by (\ref{eq:CommOrbChi}) for
$m = 2n$.  Note that our notation for the Bernoulli numbers
assumes
\begin{equation}
 \{B_1, B_2, B_3, B_4, B_5, \ldots\} = \{-\frac{1}{2}, \frac{1}{6},
              0, -\frac{1}{30}, 0, \ldots\}.
\end{equation}
This seems to be standard, but it is different from the notation
in \cite{Ahlfors}.

\subsection{The associative case}
The associative series $Q(x) = C^{3+}(x)$ is the Taylor expansion
of $f(x) = -\log(1-x) - x - x^2/2$ around $0$, which is undefined
when $x \ge 1$.  However, the Taylor expansion of $f(ix)$ is
$Q(ix)$, and $h(t,y) = \exp(-\frac1t f(it^{1/2}y))$ satisfies the
conditions of the Corollary to Theorem~\ref{thm:Exchange}. Since
the odd powers of $it^{1/2}y$ are eliminated,
Theorem~\ref{thm:Kontsevich} yields
\begin{equation}
\log \left( \denom \int_{-\infty}^\infty e^{-y^2/2}
                   e^{\frac1t f(i\sqrt{t}y)} \,dy \right)
  \approx t^{-1} \orbchi_a(-t).
\end{equation}

Now we can write
\begin{equation}
\begin{aligned}
   t^{-1}\orbchi_a(-t) &\approx
   \log\left( \frac{1}{\sqrt{2\pi}} \int_{-\infty}^{\infty}
        e^{\frac{1}{t} (-\log(1-iy\sqrt{t}) - iy\sqrt{t})} \,dy \right) \\
   &= \log\left( \frac{1}{\sqrt{2\pi}} \int_{-\infty}^{\infty}
        (1-iy\sqrt{t})^{-\frac{1}{t}} e^{-iy/\sqrt{t}} \,dy \right),
\end{aligned}
\end{equation}
where the logarithm function is defined on the complex plane minus
the positive real axis.  Next, we substitute $z =
\frac{1}{t}(iy\sqrt{t} - 1)$ to get
\begin{equation}
\begin{aligned}
   t^{-1}\orbchi_a(-t) &\approx
   \log\biggl( -\frac{i\sqrt{t}}{\sqrt{2\pi}} \int_{\Re(z) = -\frac{1}{t}}
        (-tz)^{-\frac{1}{t}} e^{-z -\frac{1}{t}} \,dz \biggr) \\
   &= \log\left( \frac{i\sqrt{t}}{\sqrt{2\pi}} (et)^{-\frac{1}{t}}
        \int_{\gamma}(-z)^{-\frac{1}{t}} e^{-z} \,dz \right),
\end{aligned}
\end{equation}
\begin{figure}[h]
\resizebox{0.3\textwidth}{!}{%
  \includegraphics{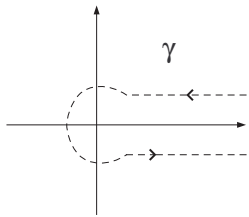}
} \caption{The contour $\gamma$}
\label{fig:HankelContour}
\end{figure}
where $\gamma$ is the curve which comes from $+\infty$ above the
positive real axis, goes around the origin counter-clockwise, and
goes back to $+\infty$ below the real axis.

This contour integral is related to the gamma function via
Hankel's formula
\begin{equation}
  \int_{\gamma} (-z)^{u-1} e^{-z} \,dz = -2i\sin(\pi u) \Gamma(u);
\end{equation}
see, for example, \cite{Whittaker}, section~12.22.  Applying this
formula, together with the identity
\begin{equation}
  \frac{\pi}{\sin{\pi u}} = \Gamma(u)\Gamma(1-u)
\end{equation}
to our integral, we get
\begin{equation}
\begin{aligned}
   t^{-1}\orbchi_a(-t) &\approx
   \log\left( \frac{2\sqrt{t}}{\sqrt{2\pi}}
       (et)^{-\frac{1}{t}} \sin\left(\pi\biggl(1-\frac{1}{t}\biggr)\right)
                                   \Gamma\biggl(1-\frac{1}{t}\biggr) \right)\\
   &= \log\left( \sqrt{2\pi t} (et)^{-\frac{1}{t}} \frac{1}{\Gamma(1/t)}
      \right).
\end{aligned}
\end{equation}
This is exactly the negative of the expression
(\ref{eq:CommOrbChi2}) we got in the commutative case!  Hence, if
we substitute $-t$ back in place of~$t$, we get the same orbifold
Euler characteristic
\begin{equation}
  \chi^{\text{orb}}_a(t) = \sum_{n=1}^{\infty}
                      \frac{B_{2n}}{2n(2n-1)} \, t^{2n}
\end{equation}
for the associative graph complex as for the commutative graph
complex.  (Which, however, is not the exact answer in the
commutative case, as we have noted earlier.)

\subsection{The Lie case}
I do not know yet how to deal with this case.  The Lie series
$Q(x) = Lie(x)$ is the Taylor expansion of $f(x) = (1-x)\log(1-x)
- 1 + x - x^2/2$ around~$0$, which is once again not defined for
$x \ge 1$, but $f(ix)$ is.  Therefore, in the same way as in the
associative case, from Theorem~\ref{thm:Kontsevich} we get
\begin{equation}
\begin{aligned}
  t^{-1}\chi^{\text{orb}}_l(-t) &\approx
  \log\left( \frac{1}{\sqrt{2\pi}}
     \int_{-\infty}^\infty
     e^{\frac1t((1-iy\sqrt{t})\log(1-iy\sqrt{t}) - 1 + iy\sqrt{t})} \,dy
                                                          \right) \\
  &= \log\left( \frac{i\sqrt{t}}{\sqrt{2\pi}}
     \int_{\Re(w)=\frac1t} e^{w\log(w)-w+w\log(t)} \,dw \right),
\end{aligned}
\end{equation}
where $w = \frac1t(1 - iy\sqrt{t})$. If we now substitute $z =
\frac{t}{e}w$, we obtain
\begin{equation} \label{eq:LieIntegral}
  t^{-1}\chi^{\text{orb}}_l(-t) \approx \log \left(
  \frac{ie}{\sqrt{2\pi t}} \int_{\Re(z)=\frac1e} z^{\frac{e}{t}z} \,dz
                                                       \right).
\end{equation}

I do not know whether this integral can be related to some known
asymptotic series, like in the other two cases; but even if it can
not, we might be able to use it to extract some information about
the coefficients of $\chi^{\textrm{orb}}_l(t)$.

\section{Summary and further questions}\label{sec:summary}
We have derived Kontsevich's formula (\ref{eq:GeneralKontsevich})
using the method of Feynman diagrams and the theory of asymptotic
expansions. Using these techniques, one could also compute the
orbifold Euler characteristic of other graph complexes, as stated
in Theorem~\ref{thm:Kontsevich}.

I have written a Mathematica notebook which performs this
computation; it can be downloaded from {\tt
http://www.math.cornell.edu/{\nohyphen}Research/{\nohyphen}Dis\-%
ser\-ta\-tions/{\nohyphen}Gerlits/{\nohyphen}code.} One enters the
formula for $Q_n/n!$ on the second line of the notebook, executes
each of the subsequent lines, and the output of the last line is
the orbifold Euler characteristic of the $Q$-graph complex for $2
\le n \le 11$ independent loops.  The numbers in
Table~\ref{table:OrbChi} were generated by this notebook.

For example, one could consider the ``chord diagram graph
complex'' by putting a chord diagram at each vertex.  Then $Q =
Ch^{3+}$, and the generating function of the orbifold Euler
characteristic is
\begin{multline*}
 -\frac{3}{8}t^2 +\frac{7}{16}t^3 -\frac{131}{128}t^4
 +\frac{449}{128}t^5
 -\frac{80179}{5120}t^6 +\frac{16459}{192}t^7 -\frac{127239605}{229376}t^8 \\
 +\frac{16956565}{4096}t^9
 -\frac{27521691751}{786432}t^{10}
 +\frac{6769184257}{20480}t^{11} + \cdots
\end{multline*}

Kontsevich claims that in the Lie case, one gets even larger
numbers for the orbifold Euler characteristic than the Bernoulli
numbers which show up in the other two cases.  This ought to be
able to be proved using (\ref{eq:LieIntegral}); also, one should
be able to give a direct proof for the fact that the orbifold
Euler characteristic of the associative graph complex is zero when
the number of loops is odd.  Another direction would be to try to
incorporate Kontsevich's orientation into the Feynman diagram
apparatus so that we could exclude graphs with
orientation-reversing automorphisms, and compute the actual
orbifold Euler characteristic in the commutative case, as well as
for graph complexes based on other species.

\bibliographystyle{plain}
\bibliography{feynman}
\end{document}